\documentclass[11pt]{article}

\usepackage{amsmath,amsfonts,amssymb}

\newcommand{\be}{\begin}
\newcommand{\Ref}[1]{(\ref{#1})}

\newcommand\fl{\mathcal R}

\newcommand{\al}{\alpha}
\newcommand{\del}{\delta}

\newcommand{\eps}{\epsilon}
\newcommand\ga{\gamma}

\newcommand\ka{\kappa}
\newcommand\lla{\lambda}
\newcommand\La{\Lambda}
\newcommand\om{\omega}

\newcommand{\si}{\sigma}
\newcommand\Si{\Sigma}

\newcommand\q{\mathbb Q}
\newcommand{\R}{\mathbb R}

\newcommand{\z}{\mathbb Z}

\newcommand{\uset}{\underset}
\newcommand{\oset}{\overset}

\newcommand{\oline}{\overline}

\newcommand{\la}{\langle}
\newcommand{\ra}{\rangle}
\newcommand{\st}{\,|\,}
\newcommand{\ti}{\tilde}

\newcommand{\prtl}{\partial}
\renewcommand\square{\kern20pt{\vbox{\hrule height.4pt
        \hbox{\vrule width.4pt height 6pt\kern6pt
                \vrule width.4pt}
        \hrule height.4pt}}}

\renewcommand\Im{\text{Im}}
\renewcommand\Re{\text{Re}}
\newcommand\rank{\text{rank}}

\newcommand\inv{^{-1}}

\newcommand\aut{\text{Aut}}

\newcommand\itr{\text{int}}

\newcommand\proof{{\em Proof.}\ }

\newcommand\SU[1]{\text{SU}(#1)}

\newcommand\SO[1]{\text{SO}(#1)}

\newcommand\Sp[1]{\text{Sp}(#1)}

\newcommand\Spin[1]{\text{Spin}(#1)}

\newcommand\vpq{V_{p,q}}
\newcommand\wpq{W_{p,q}}
\newcommand\jpq{J_{p,q}}
\newcommand\cpq{c_{p,q}}
\newcommand\fpq{f_{p,q}}

\newcommand\mij{m_{ij}}
\newcommand\mce{\mathcal E}

\begin{document}

\title{On the existence of representations of finitely presented groups
in compact Lie groups}

\author{Kim A.\ Fr\o yshov\thanks{This work was supported by the QGM
(Centre for Quantum Geometry of Moduli Spaces)
  funded by the Danish National Research Foundation.}}

\date{}

\maketitle

\be{abstract}
Given a finite, connected $2$-complex $X$ such that $b_2(X)\le1$ we
establish two existence results for representations of
the fundamental group of $X$ into
compact connected Lie groups $G$, with prescribed
values on certain loops. If $b_2(X)=1$ we assume $G=\SO3$ and that
the cup product on $H^1(X;\q)$ is non-zero.
\end{abstract}

\thispagestyle{empty}

\bibliographystyle{plain}




\section{Introduction}

In this paper we will prove two existence results for representations of
finitely presented groups into compact, connected Lie groups. The results are
most elegantly stated when the group is realized as the fundamental group of
a finite, connected $2$--dimensional CW complex $X$, as the conditions
depend only on the
cohomology of $X$. This viewpoint also reveals the relations to gauge theory
on $3$-manifolds, in which this work has its origin.

We now state our main results.
Throughout the paper we use integral coefficients for (co)homology
unless otherwise specified. The $j$'th Betti number will be denoted $b_j$.

\be{thm}\label{thm:b2-0}
Let $b_2(X)=0$, and suppose $\ga_1,\dots,\ga_r\in\pi_1(X)$ map to linearly
independent elements of $H_1(X;\q)$, where $r\ge1$.
Then for any compact, connected Lie
group $G$ and any $t_1,\dots,t_r\in G$ there exists a homomorphism
$\phi:\pi_1(X)\to G$ such that $\phi(\ga_i)=t_i$ for $i=1,\dots,r$.
\end{thm}

A reformulation of the theorem in terms of finitely presented groups is given
in Theorem~\ref{thm:finpresgr}.

Since the group homology $H_2(\pi_1(X))=0$, a theorem of Stallings
\cite[Theorem~7.4]{Stallings1} says that
the elements $\ga_1,\dots,\ga_r$ form a basis for a free subgroup of $\pi_1(X)$.
This also follows from Theorem~\ref{thm:b2-0},
because $\SO3$ contains a free subgroup
$F_2$ of rank $2$ \cite[Theorem~2.1]{wagon} and $F_2$ contains free
subgroups of any finite rank.

We now describe our second existence theorem.
Suppose $b_2(X)=1$ and let $\si$ be a generator of
$H_2(X)\approx\z$. Then
\[\mu_0:H^1(X)\times
H^1(X)\to\z,\quad(\al,\beta)\mapsto\la\al\cup\beta,\si\ra\]
is bilinear and skew-symmetric, so it defines an element
\[\mu\in\La^2(H_1(X)/T),\]
where $T\subset H_1(X)$ is the torsion subgroup.

The image of an element $\ga\in\pi_1(X)$ in $H_1(X)/T$ will be
denoted $\bar\ga$.

\be{thm}\label{thm:b2-1}
Let $b_2(X)=1$, and let $\varpi\in H^2(X;\z/2)$ with
$\la\varpi,\si\ra\neq0$. Suppose $\ga_1,\dots,\ga_{r-2}\in\pi_1(X)$
satisfy
\be{equation}\label{eq:timuga}
\mu\wedge\bar\ga_1\wedge\dots\wedge\bar\ga_{r-2}\neq0
\end{equation}
in the exterior product $\La^r(H_1(X)/T)$, where $r\ge2$.
Then there exists a homomorphism $\phi:\pi_1(X)\to\SO3$ such that
$\phi(\ga_i)=1$ for $i=1,\dots,r-2$
and $\om_2(\phi)=\varpi$. The image of such a homomorphism $\phi$
is not contained in a maximal torus in $\SO3$.
\end{thm}

The Stiefel-Whitney class $\om_2(\phi)$ is defined in Section~\ref{sec:swclass}.
Note that the Hurewicz homomorphism $\pi_2(X)\to H_2(X)$ is zero, because
$\mu\neq0$. Therefore, $H_2(\pi_1(X))\approx\z$.

Theorem~\ref{thm:b2-1} will be deduced from the more explicit
Theorem~\ref{thm:b2-1-explicit}, which
concerns representations of free groups in $\SU2$.


We now describe some relations of this paper to $3$-manifold topology.
Let $Y$ be a compact, connected, oriented $3$-manifold
with non-empty boundary. By \cite[Prop.\,4.2.7 and 4.2.13]{Gompf-Stipsicz}
$Y$ has a handle decomposition with no $3$-handles (and only one $0$-handle),
and with handles attached
in order of increasing index. Hence, $Y$ is homotopy equivalent to a
$2$-complex $X$, so Theorems~\ref{thm:b2-0} and \ref{thm:b2-1}
apply with $Y$ in place of
$X$. (In fact, the author came across these theorems while studying
moduli spaces of Bogomolny monopoles over such manifolds $Y$.)

Now suppose $b_2(Y)=1$. Then the boundary $\prtl Y$ can have at most
two components. The form $\mu_0$ can be thought of as a triple
cup product, because if $\si'\in H^1(Y,\prtl Y)$ is the Poincar\'e dual of
$\si$ (identifying the (co)homology of $X$ and $Y$) 
and $[Y]\in H_3(Y,\prtl Y)$ is the fundamental class then
\[\la\al\cup\beta,\si\ra=\la\al\cup\beta\cup\si',[Y]\ra.\]
It is not hard to see that if $\al$ lies in the image of
$H^1(Y,\prtl Y)\to H^1(Y)$ then $\iota_\al(\mu)=0$, where
$\iota_\al$ denotes contraction with $\al$. This implies that
$\mu=0$ if $\prtl Y$ is connected and $b_1(Y)<3$.

As an example, let $V$ be an oriented integral homology $3$-sphere and
$L$ an oriented link in $V$ with components $L_1,L_2$. Let $Y\subset V$ be
the complement of an open tubular neighbourhood of $L$. Then $b_2(Y)=1$, and
$\mu$ is easily computed using the formula \Ref{eqn:mu-formula} in
Section~\ref{sec:formula}. Namely,
if $a,b$ is any basis for $H_1(Y)=\z^2$ then
\[\mu=\pm\text{lk}(L_1,L_2)\,a\wedge b,\]
where $\text{lk}$ denotes the linking number. Thus, if the linking number is
non-zero then Theorem~\ref{thm:b2-1} asserts the existence of a flat
connection in the non-trivial $\SO3$ bundle over $Y$. This instance of the
theorem can also be deduced
from a result of Harper--Saveliev \cite{Harper-Saveliev2}.

It clearly suffices to prove Theorems~\ref{thm:b2-0} and \ref{thm:b2-1}
for $r=b_1(X)$.
Whereas the proof of Theorem~\ref{thm:b2-0} consists in computing
the degree of a certain diffeomorphism of $G^n$, where $n$ is the first Betti
number of the $1$-skeleton of $X$, the main point in the proof of
Theorem~\ref{thm:b2-1} is to express the quantity
\[|T|\det(\mu\wedge\bar \ga_1\wedge\cdots\wedge\bar \ga_{r-2})\]
as an intersection number in the space $Q_n^*$ of conjugacy classes of
irreducible representations $F_n\to\SU2$, where $F_n$ denotes the free
group on $n$ generators.
This is reminiscent of the definition of Casson's invariant as an intersection
number \cite{Akbulut-McCarthy}. There are also several examples of invariants
of knots or links that can be expressed as intersection numbers, see for
instance \cite{X-S-Lin,Frohman-Nicas,Harper-Saveliev1}.

The proof of Theorem~\ref{thm:b2-0} is given in Section~\ref{sec:repr-deg}.
An essential ingredient here is a result of
Gerstenhaber-Rothaus~\cite{Gerstenhaber-Rothaus} which
provides a formula for the degree of the diffeomorphism of $G^n$ defined by
an $n$-tuple of elements of $F_n$. The remaining sections of the paper are
occupied with the proof of Theorem~\ref{thm:b2-1}. In Section~\ref{sec:hom-int}
a collection of submanifolds of $Q_n^*$ is exhibited that represents a basis
for the homology group $H_{3n-6}(Q_n^*;\q)$. In
Section~\ref{sec:repr-commutators} this is used in combination with the degree
formula to compute an intersection number, thereby proving the existence of
certain representations $F_n\to\SU2$. Section~\ref{sec:swclass} describes the
second Stiefel-Whitney class of a representation $\pi_1(X)\to\SO3$ in terms
of the cellular homology of $X$. Section~\ref{sec:formula} provides a formula
for the $2$-form $\mu$ when the generator $\si$ is represented by a map from
a closed surface into $X$. In Section~\ref{sec:proof-thmb21}
these ingredients are brought together to prove
Theorem~\ref{thm:b2-1}. In the appendix a linearization map is defined
on the commutator subgroup of any group whose second rational homology
group vanishes; this may shed some light on the formula for $\mu$.

{\em Acknowledgements:} I wish to thank J\o rgen Ellegaard Andersen,
Marcel B\"okstedt, Robert Penner, Nikolai Saveliev, and Gupta Subhojoy
for valuable conversations and Takuya Sakasai for helpful
correspondence. I am also grateful to Slava Krushkal and Vladimir
Turaev for making me
aware of Stallings' theorem on free subgroups, and to
Andreas Thom for pointing out that Proposition~\ref{prop:degfw} was
proved by Gerstenhaber-Rothaus.

\section{Representations and degrees}
\label{sec:repr-deg}

Throughout this section $G$ will denote a non-trivial
compact connected Lie group
of rank $m$. 

Let $F_n$ be the free group generated by the
symbols $y_1,\dots,y_n$. We identify the abelianization of $F_n$
with $\z^n$. The image of an element $w\in F_n$ in $\z^n$ will be
denoted $\bar w$. In particular, $\bar y_i$ is the $i$'th element of the
standard basis for $\z^n$. 
Each element $w$ of $F_n$ defines a map
\[w_G:G^n\to G.\]
Namely, $w_G(g_1,\dots,g_n)$ is the image of $w$ under the unique
group homomorphism $F_n\to G$ that maps the generator $y_i$ to $g_i$.

The following proposition will be essential to the proofs
of both Theorem~\ref{thm:b2-0} and Theorem~\ref{thm:b2-1}.

\be{prop}[Gerstenhaber-Rothaus \cite{Gerstenhaber-Rothaus}]\label{prop:degfw}
If $w_1,\dots,w_n\in F_n$ then the induced map
\[f:=(w_1)_G\times\cdots\times(w_n)_G:G^n\to G^n\]
has degree $\deg(f)=(\det(\bar w_1,\dots,\bar w_n))^m$.\square
\end{prop}

The case $n=1$ is a classical theorem of Hopf~\cite{Hopf2}
(see also \cite[p.174]{GHV2}). The general case
can be proved by combining Hopf's theorem with
the fact (see \cite[p.\,169]{GHV2}) that the two maps
$\mu,\mu':G\times G\to G$ given by
\[\mu(g,h)=gh,\quad\mu'(g,h)=hg\]
induce the same map in rational cohomology. Alternatively -- and this is
the approach taken in \cite{Gerstenhaber-Rothaus} -- one
can compute the effect of $f$ on the top cohomology of $G^n$ by means of
another theorem of Hopf~\cite{Hopf1} which says that the rational
cohomology of $G$
is an exterior algebra generated by homogeneous elements of odd degree.

\be{thm}\label{thm:finpresgr}
Let $H$ be a
group with generators $y_1,\dots,y_n$ and relations $w_1,\dots,w_s$,
where $\bar w_1,\dots,\bar w_s$ are linearly independent. Let $A:=H/[H,H]$
be the abelianization of $H$.
If $h_1,\dots,h_r\in H$ map to linearly independent elements of
$A\otimes\q$ then for
any $t_1,\dots,t_r\in G$ there exists a homomorphism $\phi:H\to G$
such that $\phi(h_i)=t_i$ for $i=1,\dots,r$.
\end{thm}

\proof We may assume $r=\rank(A)$.
Choose a lift $z_i\in F_n$ of $h_i$ for each $i$ and let
\[f:G^n\to G^n\]
be the map defined by the $n$ elements
$w_1,\dots,w_s,z_1,\dots,z_r\in F_n$ as in Proposition~\ref{prop:degfw}.
Then the set of homomorphisms $\phi$ as in the theorem can be
identified with $f\inv(1,\dots,1,t_1,\dots,t_r)$. Thus, it suffices to
show that $f$ has non-zero degree, since this implies that $f$ is
surjective.

Let $N$ be the kernel of the projection $\pi:F_n\to H$, and let $K$ be
the kernel of the induced homomorphism $\tau:\z^n\to A$. Then we have a
commutative diagram
\[\be{array}{ccccccccc}
1 & \longrightarrow & N & \longrightarrow & F_n &
\oset\pi\longrightarrow & H & \longrightarrow & 1\\
& & \si'\downarrow\hphantom{\si'} & &
\si\downarrow\hphantom{\si} & &
\si''\downarrow\hphantom{\si''} & &\\
0 & \longrightarrow & K & \longrightarrow & \z^n &
\oset\tau\longrightarrow & A & \longrightarrow & 0
\end{array}\]
where $\si,\si''$ are the abelianization maps. Because $\pi$ maps the
commutator subgroup of $F_n$ onto the commutator subgroup of $H$, the
map $\si'$ is surjective. But this just means that $K$ is generated by
$\bar w_1,\dots,\bar w_s$, hence these elements form a basis for $K$.

By the elementary divisors theorem \cite{Lang2} there is a basis $b_1,\dots,b_n$
for $\z^n$ and integers $p_1,\dots,p_s$ such that $p_1b_1,\dots,p_sb_s$ is a
basis for $K$. This implies that $b_{s+1},\dots,b_n$ map to a basis
for $A/T$, where $T$ is the torsion subgroup of $A$, and
\[\bar w_1\wedge\cdots\wedge\bar w_s=\pm p_1\cdots
p_sb_1\wedge\cdots\wedge b_n=\pm|T|b_1\wedge\cdots\wedge b_n.\]
Therefore, by Proposition~\ref{prop:degfw} we have
\be{align*}
\deg(f)&=(\det(\bar w_1,\dots,\bar w_s,\bar z_1,\dots,\bar z_r))^m\\
&=\pm(|T|\det(b_1,\dots,b_s,\bar z_1,\dots,\bar z_r))^m\\
&=\pm(|T|\det(\bar h_1,\dots\bar h_r))^m\\
&\neq0,
\end{align*}
where $\bar h$ denotes the image in $A/T$
of an element $h\in H$.\square

{\em Proof of Theorem~\ref{thm:b2-0}:} Let
\[C_2\oset\prtl\to C_1\oset\prtl\to C_0\]
be the cellular chain complex of $X$, so that $C_k=H_k(X^k,X^{k-1})$,
where $X^k$ is the $k$-skeleton of $X$. Then $X^1$ is homotopy
equivalent to a wedge sum $\vee_nS^1$ of $n$ circles for some
$n\ge0$. Let $x_0\in X^1$ and $d_0\in S^1$ be base-points and fix a basis
$y_1,\dots,y_n$ for the free group $\pi_1(X^1,x_0)$.

Let the $2$-cells in $X$ be numbered from $1$ to $s$. For
$i=1,\dots,s$ choose a base-point preserving map $\ell_i:S^1\to X^1$ that 
is homotopic to the attaching map of the $i$'th $2$-cell. Let
$w_i\in \pi_1(X^1,x_0)=F_n$ be the class represented by $\ell_i$. Then
$H:=\pi_1(X,x_0)$ has the presentation
\be{equation*}
H=\la y_1,\dots,y_n\st w_1,\dots,w_s\ra.
\end{equation*}
Moreover, if $a_i\in C_2$ is the generator corresponding to the $i$'th
$2$-cell then
\[\prtl a_i=\bar w_i\in\z^n=H_1(X^1),\]
where we identify $H_1(X^1)$ with its image in $H_1(X^1,X^0)$. Thus,
the assumption $b_2(X)=0$ means precisely that $\bar w_1,\dots,\bar
w_s$ are linearly independent, so the theorem follows from
Theorem~\ref{thm:finpresgr}.\square

\section{Homology and intersection numbers}
\label{sec:hom-int}

This section is concerned with the homology of certain quotient spaces $Q_n^*$
associated to the Lie group $\Sp1$ of unit quaternions. (Of course,
$\Sp1\approx\SU2$.)

We begin by explaining our orientation conventions. If $W$ is a smooth oriented
manifold and $Z\subset W$ an oriented submanifold then the normal
bundle $NZ$ will be oriented so that for any Euclidean metric in the
tangent bundle $TW$ the isomorphism
\[TW|_Z=NZ\oplus TZ\]
preserves orientations. If $V$ is another smooth oriented manifold and
$f:V\to W$ is transverse to $Z$ then the submanifold $f\inv Z$ will be
oriented such that the isomorphism $N(f\inv Z)\to f^*(NZ)$ preserves
orientations. Given a smooth fibre bundle $\pi:P\to V$ where $P$ is
oriented as manifold, the fibres $F$ of $\pi$ will be oriented
such that for any Euclidean metric in $TP$ the isomorphism
\[TP|_F=TF\oplus\pi^*(TV)|_F\]
preserves orientations. The latter convention prevents some
annoying signs. In particular, if $p\in F$ is a regular value of
$f:V\to F$ and $s$ the corresponding section of the product bundle
$F\times V\to V$ then $f\inv(p)=s\inv(\{p\}\times V)$ as oriented manifolds.

We shall need an extension of the notion of intersection number defined in
\cite{Guillemin-Pollack}. Let $V,W,Z$ be as above. Suppose $Z$ is
closed as a subset of $W$ and that
\[\dim\,V+\dim\,Z=\dim\,W.\]
For any smooth map $f:V\to W$ such that $f\inv Z$ is compact
we define the intersection number $I(f,Z)$
as follows. Let $g:V\to W$ be any smooth map transverse to $Z$ such that
$g$ is {\em compactly homotopic to $f$}, by which we mean that $g$ is
homotopic to $f$ relative to the complement of a compact subset of $V$.
Then we set
\[I(f,Z):=\#g\inv Z,\]
where $\#$ denotes the number of points counted with sign.
If $V$ is in fact a submanifold of $W$ and $f$ is the inclusion map
then we define the intersection number
$V\cdot Z:=I(f,Z)$.

Throughout this section we will use the notation
\be{equation}\label{eqn:GG}
G:=\Sp1,\quad G':=\Sp1/\pm1=\SO3.
\end{equation}

For any integer $n\ge2$ we define a right $G'$--space $P_n$ as follows.
Let $P_n:=G^n$ as smooth, oriented manifold, and let $g\in G$ act on $P_n$ by
conjugation with $g\inv$ in each factor. Since $-1\in G$ acts
trivially, this right action of $G$ descends to a right action of $G'$
on $P_n$. Let $P_n^*$ be the open subset of $P_n$ consisting of those
points in $P_n$ that have trivial stabilizer in $G'$. Since $P_n$ is
Hausdorff and $G'$ is compact, the quotient space
\[Q_n:=P_n/G'\]
is Hausdorff. Let $Q_n^*$ be the image of $P_n^*$ in $Q_n$, which is
an open subset of $Q_n$. Then $Q_n^*$ has a canonical smooth structure
such that $\pi:P_n^*\to Q_n^*$ is a principal
$G'$--bundle. It is easy to see that $Q_2^*$ is diffeomorphic to
$\R^3$. Each fibre of $\pi$ inherits an orientation from $G'$, and
since $G'$ is connected, the base manifold $Q_n^*$ is orientable. Let
$Q_n^*$ have the orientation compatible with the orientations of the
total space and fibres of $\pi$ as stipulated above.

\be{prop}\label{prop:bq}
$b_{3n-6}(Q_n^*)=n(n-1)/2$.
\end{prop}

\proof In this proof we use real coefficients for (co)homology.
The proposition clearly holds for $n=2$, so let $n\ge3$. 
Set
\[R_n:=P_n\setminus P_n^*,\]
so that $R_n$ consists of those $n$-tuples of elements from $G$ such
that all elements are contained in the same maximal torus.
We will show that
\[b_{3n-6}(Q_n^*)=b_{3n-3}(P_n^*)=b_2(R_n)=n(n-1)/2.\]
The last equality is a special case of results by Baird \cite{Baird1}, but
we include a direct proof here for completeness.
Identifying $S^2$
with the unit sphere in the space of pure quaternions (ie
quaternions with zero real part) we have a surjective map
\[\chi_n:S^2\times T^n\to R_n,\quad(x,z)\mapsto(\Re(z_j)+\Im(z_j)\cdot
x)_{j=1,\dots,n},\]
where $T^n:=U(1)^n$ and $z=(z_1,\dots,z_n)$. 
Clearly, $\chi_n$ factors through a map $\ti R_n\to R_n$, where
\[\ti R_n:=(S^2\times T^n)/(x,z)\sim(-x,z\inv),\]
which is a smooth manifold.
For each $\eps\in\{-1,1\}^n$ let $B_\eps\approx\mathbb{RP}^2$ be the
image of $S^2\times\{\eps\}\subset S^2\times T^n$ in $R^n$. Then $R_n$
is obtained from $\ti R_n$ by collapsing each $B_\eps$ to
a point $x_\eps$, where $x_\eps\neq x_{\eps'}$ if $\eps\neq\eps'$. Let $N$ be
a closed tubular neighbourhood of $B:=\cup_\eps B_\eps$ in $\ti R_n$ and set
$Z:=\ti R_n\setminus\itr(N)$. Because each $B_\eps$ is a rational homology
ball it follows from the Mayer-Vietoris sequences for $(N,Z)$ and the image of
$(N,Z)$ in $R_n$ that the projection $\ti R_n\to R_n$ induces an
isomorphism
\[H_k(\ti R_n)\oset\approx\to H_k(R_n)\]
for each $k$.
Now, $H^*(\ti R_n)$ is isomorphic to the $1$-eigenspace of the endomorphism of
$H^*(S^2\times T^n)$ induced by the involution $(x,z)\to(-x,z\inv)$.
This yields in particular
\[b_2(R_n)=n(n-1)/2.\]
Using the fact that $\chi_1:S^2\times S^1\to R_1=S^3$ has degree~$2$, one finds
that the inclusion $R_n\to P_n$ induces a surjection $H_3(R_n)\to H_3(P_n)$.
Since $H_2(P_n)=0$ the homology sequence of the pair $(P_n,R_n)$ shows that the
connecting homomorphism
\[H_3(P_n,R_n)\oset\approx\to H_2(R_n)\]
is an isomorphism.

Now observe that $P_n$ and $R_n$ are real algebraic subsets of $\R^{4n}$. This is
obvious for $P_n$. To see that it holds for $R_n$, note that
\[R_n=\{(g_1,\dots,g_n)\in P_n\st\text{$g_pg_q=g_qg_p$ for $1\le p<q\le n$.}\}\]
Therefore, by \cite{Hironaka} both $P_n$ and $R_n$
admit triangulations. (In the case
of $P_n$ this is of course elementary.) By
\cite[Cor.\,6.1.11 and Thm.\,6.2.17]{Spanier} we have a Poincar\'e duality
isomorphism
\[H^3(P_n,R_n)\approx H_{3n-3}(P_n^*).\]

We now apply the homology spectral sequence of the fibration $P_n^*\to Q_n^*$.
Since $P_n$ is simply-connected and
\[\dim\,P_n-\dim(S^2\times T^n)=2n-2\ge4,\]
one can show by a transversality argument that
$P_n^*=P_n\setminus\chi_n(S^2\times T^n)$ is simply-connected.
Because the fibre $G'$ of the bundle $P_n^*\to Q_n^*$ is path-connected, we
conclude that $Q_n^*$ is simply-connected. Hence the fibration is orientable
and the $E^2$-page of the spectral sequence is
\[E^2_{p,q}\approx H_p(Q_n^*;H_q(G')).\]
Since $Q_n^*$ is a non-compact $(3n-3)$-manifold, we have $H_{3n-3}(Q_n^*)=0$.
Therefore,
\[H_{3n-3}(P_n^*)\approx E^\infty_{3n-6,3}\approx E^2_{3n-6,3}\approx H_{3n-6}(Q_n^*)\]
and the proposition is proved.\square

For each pair $p,q$ of integers satisfying $1\le p<q\le n$ we have a
$G'$--equivariant embedding
\be{equation}\label{eqn:q2qn}
P_2\to P_n,\quad(g,h)\mapsto(1,\dots,1,g,1,\dots,1,h,1,\dots,1),
\end{equation}
where $g$ and $h$ appear in the $p$'th and $q$'th place, resp.
This map induces a topological embedding
$I_{p,q}:Q_2\to Q_n$
that restricts to a smooth embedding $Q_2^*\to Q_n^*$
whose image we denote by $\wpq$. Clearly, $\wpq=I_{p,q}(Q_2)\cap
Q_n^*$ is a closed subset of $Q_n^*$.
Letting $i,j$ denote the usual
anti-commuting quaternions for the time being, the composition of the map
\be{equation}\label{eqn:jpqdef}
\be{aligned}
\jpq:P_{n-2}&\to P_n^*,\\
(g_1,\dots,g_{n-2})&\mapsto
(g_1,\dots,g_{p-1},i,g_p,\dots,g_{q-2},j,g_{q-1},\dots,g_{n-2}).
\end{aligned}
\end{equation}
with the projection $P_n^*\to Q_n^*$ is a smooth embedding $P_{n-2}\to
Q_n^*$ whose image we denote by $\vpq'$.
Let $\wpq$ and $\vpq'$ have the orientations
inherited from $Q_2^*$ and $P_{n-2}$, resp.

For each positive integer $m$ let $G$ act on $G^m$ by conjugation in
each factor. This (left) action descends to an action of $G'$ on $G^m$ and we
have an associated bundle of Lie groups
\[E_{m,n}:=P_n^*\uset{G'}\times G^m\to Q_n^*.\]
Because the action of $G'$ on $G^m$ preserves orientations, each fibre
of $E_{m,n}$ inherits an orientation from $G^m$. Since the base space
$Q_n^*$ is oriented, we get an induced orientation on the total space
$E_{m,n}$ by the above convention.
The fixed-point $-1\in G$ gives rise to a section of
$E_{1,n}$ whose image we denote by $-\mathbf1$. 

For $p,q$ as above the commutator map
\[P_n\to G,\quad(g_1,\dots,g_n)\mapsto[g_p,g_q]\]
is $G'$--equivariant (in the sense that it intertwines the right action
of $h\inv$ on $P_n^*$ with the left action of $h$ on $G$, for $h\in
G'$) and therefore defines a section $\zeta_{p,q}$ of $E_{1,n}$. At
this point we recall the following well known facts about $G$:

\be{lemma}\label{lemma:quat}
\be{description}
\item[(i)]Two unit quaternions anti-commute if and only if their real
  parts vanish and their imaginary parts are orthogonal.
\item[(ii)]$1\in G$ is the only singular value of the commutator map
\[P_2\to G,\quad(g,h)\mapsto[g,h].\square\]
\end{description}
\end{lemma}

It follows from the second part of the lemma that $\zeta_{p,q}$ is
transverse to $-\mathbf1$. Set
\[\vpq:=\zeta_{p,q}\inv(-\mathbf1)\subset Q_n^*\]
as oriented manifold. 

\be{lemma}\label{lemma:V122}
$V_{1,2}\subset Q_2^*$ consists of a single point, which is positively oriented.
\end{lemma}

\proof From Lemma~\ref{lemma:quat}~(i) we see that $V_{1,2}$
consists of a single point. The sign 
can be determined by an explicit computation, which we omit.\square

\be{prop}\label{prop:wvint}
In $Q_n^*$ the following hold:
\be{description}
\item[(i)]The submanifolds $\wpq$ and $\vpq$
intersect transversely in a single point, which is positively oriented.
\item[(ii)]$\wpq$ and $V_{p',q'}$ are disjoint for $(p,q)\neq(p',q')$.
\item[(iii)]$\vpq'=(-1)^{p+q+1}\vpq$.
\end{description}
\end{prop}

\proof Part~(i) is a consequence of Lemma~\ref{lemma:V122}, whereas
(ii) is trivial. Lemma~\ref{lemma:quat}~(i) implies that $\vpq'=\vpq$ as
smooth manifolds.
On the other hand, one easily checks that
$\wpq\cdot\vpq'=(-1)^{p+q+1}$, proving (iii).
\square

It follows from the proposition that the $\vpq$'s represent linearly
independent classes in $H_{3n-6}(Q_n^*;\q)$. Combining this with
Proposition~\ref{prop:bq} we obtain

\be{cor}\label{cor:basis}
The oriented submanifolds $\{\vpq\}_{1\le p<q\le n}$ of $Q_n^*$
represent a basis for
$H_{3n-6}(Q_n^*;\q)$.
\end{cor}

Now let $n=2\rho$, where $\rho$ is a positive integer. The
$G'$--equivariant map
\[c_\rho:P_{2\rho}\to G,\quad(g_1,\dots,g_{2\rho})\mapsto
\prod_{\ell=1}^\rho[g_{\ell},g_{\ell+\rho}]\]
defines a section $s_\rho$ of $E_{1,2\rho}$, which is transverse to
$-\mathbf1$ by Lemma~\ref{lemma:quat}~(ii). Set
\[M_\rho:=s_\rho\inv(-\mathbf1)\subset Q_{2\rho}^*\]
as oriented manifold. Because $c_\rho\inv(-1)$ is a compact subset of
$P_n^*$, the space $M_\rho$ is compact. 

The spaces $M_\rho$ have been the subject of much research,
see for instance \cite{Thaddeus1} and the references therein.
The following basic result might be known, but
we have been unable to find a reference.

\be{prop}\label{prop:Mclass}
The class in $H_{6\rho-6}(Q_{2\rho}^*;\q)$ represented by
$M_\rho$ is given by
\[[M_\rho]=\sum_{\ell=1}^\rho[V_{\ell,\ell+\rho}].\]
\end{prop}

\proof It follows from Lemma~\ref{lemma:V122} that
\[W_{\ell,\ell+\rho}\cdot M_\rho=1,\]
and one clearly has $\wpq\cap M_\rho=\emptyset$ unless
$(p,q)=(\ell,\ell+\rho)$ for some $\ell$. The proposition now follows
from Corollary~\ref{cor:basis} and Proposition~\ref{prop:wvint}.\square

\section{Representations and commutators}
\label{sec:repr-commutators}

We will now use the results of Section~\ref{sec:hom-int} to prove an existence
theorem for representations of free groups into $\Sp1$, from which we will
deduce Theorem~\ref{thm:b2-1}.

Suppose $v_0,\dots,v_k$ are elements of the free group
$F_n$ such that $v_0$ is a
product of commutators,
\be{equation}\label{eqn:v0-comm}
v_0=\prod_{\ell=1}^\rho[u_\ell,u_{\ell+\rho}],
\end{equation}
where $u_1,\dots,u_{2\rho}\in F_n$. 
Set $L:=\z^n$ and
\[\lla:=\sum_{\ell=1}^\rho\bar u_\ell\wedge\bar u_{\ell+\rho}\in\La^2L.\]
In the appendix it is shown that $\lla$ is in fact determined by $v_0$, but
we will not need this fact.
Let $\eps_0,\dots,\eps_k\in\{\pm1\}$ with $k\ge0$ and $\eps_0=-1$.

\be{thm}\label{thm:b2-1-explicit}
Suppose
\[\lla\wedge\bar v_1\wedge\cdots\wedge\bar v_k\neq0\quad\text{in
  $\La^{k+2}L$}.\]
Then there exists a homomorphism $\psi:F_n\to\Sp1$ satisfying
\[\psi(v_i)=\eps_i,\quad i=0,\dots,k.\]
\end{thm}

As an application, let $H$ be the
group with generators $y_1,\dots,y_n$ and relations
$v_0,\dots,v_k$. Then $\psi$ induces a homomorphism
$\phi:H\to\SO3$ such that
\be{itemize}
\item the image of $\phi$ is not contained in a maximal torus of
  $\SO3$,
\item $\phi$ does not lift to a homomorphism $H\to\Sp1$.
\end{itemize}

{\em Proof of Theorem~\ref{thm:b2-1-explicit}:} 
We may assume $k=n-2$.
The set of equivalence classes of such
representations $\psi$ form a subset $\fl\subset Q_n$. In fact, because $v_0$
is a product of commutators we have $\fl\subset Q^*_n$.
We will express $\fl$ in 
a different way, making use of \Ref{eqn:v0-comm}. Let
$f_0:Q_n\to Q_{2\rho}$ 
be the map defined by $u_1,\dots,u_{2\rho}$.
Then $U:=f_0\inv Q_{2\rho}^*$ is an open subset of $Q_n^*$. Then
\[f:=f_0|_U:U\to Q_{2\rho}^*\]
is a proper map.
The elements $(v_1,\dots,v_{n-2})$ define a $G'$-equivariant map
\[h:P_n\to G^{n-2}\]
which in turn determines a section $\xi_0$ of $E_{n-2,n}$.
Set $\xi:=\xi_0|_U$. Let $\mce$ denote the image of the section of $E_{n-2,n}$
corresponding to the fixed-point $(\eps_1,\dots,\eps_{n-2})$ of $G^{n-2}$.
Then
\[\fl=(f\times\xi)\inv(M_\rho\times\mce).\]
We will show that the intersection number
\[\ka:=I(f\times\xi,M_\rho\times\mce)\]
is non-zero, which will imply that $\fl$ is non-empty. Note that the
intersection number is indeed well-defined, since $M_\rho$ and $\mce$ are closed
subsets of $Q_{2\rho}^*$ and $E_{n-2,n}$ resp.\ and 
$f\inv M_\rho$ is a compact subset of $U$.

\be{prop}\label{prop:ka}
$\ka=\det(\lla\wedge\bar v_1\wedge\cdots\wedge\bar v_{n-2})\neq0$.
\end{prop}

Here, `$\det$' denotes the standard isomorphism $\La^nL\to\z$.

\proof Choose a smooth map $f_1:U\to Q_{2\rho}^*$ that is transverse to $M_\rho$
and compactly homotopic to $f$. (See the beginning of Section~\ref{sec:hom-int}
for the definition of `compactly homotopic'.) Then
$N:=f_1\inv M_\rho$ is a compact, oriented $(3n-6)$-dimensional submanifold of
$U$. Choose a smooth section $\xi_1$ of $E_{n-2,n}|_U$ compactly homotopic to
$\xi$ such that $\xi_1|_N$ is transverse to $\mce$.
Then $f_1\times\xi_1$ is transverse to $M_\rho\times\mce$, so
\be{align*}
\ka&=\#(f_1\times\xi_1)\inv(M_\rho\times\mce)\\
&=\#(\xi_1|_N)\inv\mce\\
&=I(\xi|_N,\mce).
\end{align*}

For any $v\in F_n$ and $1\le p\le n$ let $\deg_p(v)$ denote the $p$'th
component of $\bar v\in\z^n$.

\be{lemma}
The class in $H_{3n-6}(Q_n^*;\q)$ represented by $N$ is given by
\be{equation}\label{eqn:cpq}
[N]=\sum_{p<q}\cpq[\vpq],
\end{equation}
where
\[\cpq=\sum_{\ell=1}^\rho[\deg_p(u_\ell)\deg_q(u_{\ell+\rho})
-\deg_q(u_\ell)\deg_p(u_{\ell+\rho})].\]
\end{lemma}

\proof That $[N]$ can be expressed in the form \Ref{eqn:cpq} follows from
Corollary~\ref{cor:basis}. Set $\wpq^-:=\wpq\cap U$. Then
$\fpq:=f|_{\wpq^-}:\wpq^-\to Q_{2\rho}^*$ is a proper map and
\[\cpq=\wpq\cdot N=\wpq^-\cdot N=I(\fpq,M_\rho)=
\sum_{\ell=1}^\rho I(\fpq,V_{\ell,\ell+\rho}).\]

Let $\psi_{p,q,\ell}:P_2\to P_2$ be the composition of the embedding
$P_2\to P_n$ in 
\Ref{eqn:q2qn} with the map $P_n\to P_2$ defined by $(u_\ell,u_{\ell+\rho})$.
From Lemma~\ref{lemma:V122} and Proposition~\ref{prop:degfw} we obtain
\be{align*}
I(\fpq,V_{\ell,\ell+\rho})&=\deg(\psi_{p,q,\ell})\\
&=\deg_p(u_\ell)\deg_q(u_{\ell+\rho})
-\deg_q(u_\ell)\deg_p(u_{\ell+\rho}).
\end{align*}
To prove the first of the above two equalities we choose a smooth map
$b:\wpq^-\to Q_{2\rho}^*$ transverse to $V_{\ell,\ell+\rho}$ and compactly
homotopic to $\fpq$. We then observe that any such homotopy can be
lifted to a $G'$-equivariant smooth map
$B:P_n^*|_{\wpq^-}\times[0,1]\to P_{2\rho}^*$ such that $B(\,\cdot\,,0)$ is
the map induced by $u_1,\dots,u_{2\rho}$.\square

Returning to the proof of Proposition~\ref{prop:ka} we have
\be{align*}
\ka&=I(\xi|_N,\mce)=\sum_{p<q}\cpq I(\xi|_{\vpq},\mce)\\
&=\sum_{p<q}(-1)^{p+q+1}\cpq I(\xi|_{\vpq'},\mce),
\end{align*}
where the second equality comes from Proposition~\ref{prop:wvint}\,(iii).
Set
\[V:=\bar v_1\wedge\dots\wedge\bar v_{n-2}.\] 
Because $E_{n-2,n}|_{\vpq'}$ is trivial we have
\[I(\xi|_{\vpq'},\mce)=\deg(h\circ\jpq)
=(-1)^{p+q+1}\det(e_p\wedge e_q\wedge V),\]
where $\jpq:P_{n-2}\to P_n^*$ is the map in \Ref{eqn:jpqdef} and $e_1,\dots,e_n$
is the standard basis for $L=\z^n$. Thus,
\[\ka=\sum_{p<q}\cpq\det(e_p\wedge e_q\wedge V)
=\det(\lla\wedge V).\square\]

\section{The second Stiefel--Whitney class of an $\SO3$
representation}
\label{sec:swclass}

Let $X$ be a finite, connected CW-complex with base-point $x_0$, and let
\[\phi:H:=\pi_1(X,x_0)\to\SO3\]
be any homomorphism. 
Let $\ti X\to X$ be the universal covering. A choice of base-point in $\ti X$
lying above $x_0$ makes $\ti X$ a principal $H$-bundle, so
we can associate to $\phi$ an $\SO3$-bundle $E_\phi\to X$, whose
second Stiefel-Whitney class we denote by $\om_2(\phi)\in H^2(X;\z/2)$. The aim
of this section is to describe $\om_2(\phi)$ in terms of its value on
any element of
$H_2(X;\z/2)$. (A description of $\om_2(\phi)$ in terms of group
cohomology can be found in \cite[Section~3.1]{rub-sav0}.)

We will use the same notation as in the proof of
Theorem~\ref{thm:b2-0} in Section~\ref{sec:repr-deg},
so in particular, $C_*$ denotes the
cellular chain complex of $X$. Let 
$c=\sum_ia_i\otimes c_i$ be any cycle in $C_2\otimes\z/2$, where $c_i\in\z/2$.
For each $j$ choose a lift 
$q_j\in\Sp1$ of $\phi(y_j)\in\SO3$, and let $\psi:F_n\to\Sp1$ be the
unique homomorphism such that $\psi(y_j)=q_j$. Then $\psi(w_i)=\pm1$
for each $i$.

\be{prop}\label{prop:w2descr}
Let $\del$ denote the group isomorphism $\{\pm1\}\to\z/2$. Then
\[\la \om_2(\phi),[c]\ra=\sum_{i=1}^sc_i\del(\psi(w_i)).\]
\end{prop}

Here $[c]\in H_2(X;\z/2)$ denotes the homology class of $c$.
The proposition holds more generally for $\SO N$ representations,
$N\ge3$, if one replaces $\Sp1$ by $\Spin N$.

\proof Let $J$ be the set of indices $i\in\{1,\dots,s\}$ such that $c_i=1$,
and let $d_0\in S^1=\prtl D^2$ be a base-point as before.
Choose an embedding $h:D^2\times J\to\itr(D^2)$ and a map
$g:D^2\to X$ such that $g(d_0)=x_0$ and the following properties hold:
\be{itemize}
\item $g$ maps the complement of the image of $h$ into $X^1$,
\item for each $i\in J$ the map $g(h(\cdot,i)):D^2\to X$ is the
  characteristic map of the $i$'th $2$-cell of $X$,
\item $g|_{S^1}$ represents the element $\prod_{i\in J}w_i$ of $\pi_1(X^1,x_0)$,
where the product is taken according to some ordering of $J$.
\end{itemize}

Let $Y$ be the result of attaching a $2$-cell to $X^1$ with $g|_{S^1}$ as
attaching map. Let $f:Y\to X$ be the cellular map obtained by combining
$g$ with the inclusion $X^1\to X$. Then $f_*:H_2(Y;\z/2)\to H_2(X;\z/2)$
maps the non-zero element of $H_2(Y;\z/2)$ to $[c]$.
Set $K:=\pi_1(Y,x_0)$ and let
$\ga$ denote the composite homomorphism
\[K\oset{f_*}\to H\oset\phi\to\SO3.\]
Let $\ti Y\to Y$ be the universal covering, and choose a base-point in $\ti Y$
lying above $x_0$.
The base-point preserving map $\ti Y\to\ti X$ covering $f$ is equivariant
with respect to $f_*:K\to H$ and induces an isomorphism
$E_\ga\to f^*E_\phi$.
It is now easy to verify the following equivalences:
\be{align*}
\la \om_2(\phi),[c]\ra=0 & \iff \om_2(\ga)=0\\
& \iff\text{$E_\ga$ lifts to an $\Sp1$ bundle}\\
& \iff\text{$\ga$ lifts to a homomorphism $K\to\Sp1$}\\
& \iff\psi(\textstyle\prod_{i\in J}w_i)=1\\
& \iff\sum_{i=1}^sc_i\del(\psi(w_i))=0,
\end{align*}
where in the fourth equivalence we used the fact that $\psi(\prod_{i\in J}w_i)$
is independent of the choice of the lifts $q_j$ since $c$ is a cycle.\square

\be{prop}\label{prop:max-torus}
If the image of $\phi$ is contained in a maximal torus of $\SO3$, then
$\la \om_2(\phi),z\ra=0$ for every $z\in H_2(X;\z)$.
\end{prop}

\proof Let $z$ be represented by a cycle $c=\sum_{i=1}^sc_ia_i$ in
$C_2$, where $c_i\in\z$.
Then $w:=\prod_i(w_i)^{c_i}\in F_n$ is a product of commutators. Since $\psi$
takes values in a maximal torus of $\Sp1$, we have $\psi(w)=1$, so
\[0=\del(\psi(w))=\sum_ic_i\del(\psi(w_i))=\la \om_2(\phi),z\ra.
\square\]

\section{A formula for $\mu$}
\label{sec:formula}

We now give a formula for the $2$-form $\mu$ in Theorem~\ref{thm:b2-1}
that is useful for computations. Suppose
the generator $\si\in H_2(X)$ is represented by a map $f:\Si\to X$,
where $\Si$ is a closed, connected surface of genus $\rho$. Let
$b_1,\dots,b_{2\rho}$ be a symplectic basis for
$H_1(\Si)$, so for $1\le i<j\le2\rho$ one has the intersection numbers
\[b_i\cdot b_j=\be{cases}
1 & \text{if $j=i+\rho$,}\\
0 & \text{else.}
\end{cases}\]
For $\al,\beta\in H^1(\Si)$ one easily checks that
\[\la\al\cup\beta,[\Si]\ra=
\sum_{\ell=1}^\rho\left[\la\al,b_\ell\ra\la\beta,b_{\ell+\rho}\ra
-\la\beta,b_\ell\ra\la\al,b_{\ell+\rho}\ra\right].\]
Letting $\ti b_j$ denote the image of $f_*(b_j)$ in $H_1(X)/T$ we obtain
\[\la\al\cup\beta,\si\ra=
\sum_{\ell=1}^\rho\left[\la\al,\ti b_\ell\ra\la\beta,\ti b_{\ell+\rho}\ra
-\la\beta,\ti b_\ell\ra\la\al,\ti b_{\ell+\rho}\ra\right],\]
so that
\be{equation}\label{eqn:mu-formula}
\mu=\sum_{\ell=1}^\rho\ti b_\ell\wedge\ti b_{\ell+\rho}.
\end{equation}

\section{Proof of Theorem~\ref{thm:b2-1}}
\label{sec:proof-thmb21}

We may assume $r=b_1(X)$.
If $\phi$ is a homomorphism as in the theorem then by
Proposition~\ref{prop:max-torus} the image of $\phi$ cannot be
contained in a maximal torus. We now prove the existence of $\phi$.
We continue using the notation introduced in the proof of
Theorem~\ref{thm:b2-0}, except that we now assume the $2$-cells of $X$ are
numbered from $0$ to $s$, with corresponding generators $a_0,\dots,a_s$
of $C_2$, so that $H$ has the presentation
\[H=\la y_1,\dots,y_n\st w_0,\dots,w_s\ra.\] 

We begin by modifying the relations in this presentation.
Choose an invertible integral matrix
$M=(\mij)_{0\le i,j\le s}\in\text{GL}_{s+1}(\z)$ with
\[\sum_{i=0}^sm_{i0}a_i=\si.\]
Choose
an automorphism $\al\in\aut(F_{s+1})$ that maps to $M$ under the canonical
surjective homomorphism $\aut(F_{s+1})\to\text{GL}_{s+1}(\z)$. Let $p:F_{s+1}\to
F_n$ be the unique homomorphism with $p(y_{i+1})=w_i$ for $i=0,\dots,s$,
and set
\be{equation}\label{eqn:wipdef}
w_i':=p(\al(y_{i+1})),\quad i=0,\dots,s.
\end{equation}
Thus, each $w_i'$ can be expressed as a word in $w_0,\dots,w_s$,
and conversely, each $w_i$ is a word in $w_0',\dots,w_s'$.
Hence, we have the new presentation
\[H=\la y_1,\dots,y_n\st w_0',\dots,w_s'\ra.\]

Choose $\eta_0,\dots,\eta_s\in\z/2$ such that for any cycle
$c=\sum_{i=0}^sa_i\otimes c_i$ in $C_2\otimes\z/2$ one has
\[\la\varpi,c\ra=\sum_ic_i\eta_i.\]
For $j=0,\dots,s$ set
\[\eps_j:=\del\inv(\sum_i\mij\eta_i)\in\{\pm1\},\]
where $\del$ is as in Proposition~\ref{prop:w2descr}. Then
\[\eps_0=\del\inv\la\varpi,\si\ra=-1.\]

\be{lemma}\label{lemma:tphi-eps}
If $\psi:F_n\to\Sp1$ is any homomorphism such that $\psi(w_j')=\eps_j$ for
each $j$ then the induced homomorphism $\phi:H\to\SO3$ satisfies
$\om_2(\phi)=\varpi$.
\end{lemma}

\proof For each $j$ one has
\[\sum_i\mij\eta_i=\del(\psi(w_j')=\sum_i\mij\del(\psi(w_i)),\]
hence $\eta_i=\del(\psi(w_i))$ and the lemma follows from
Proposition~\ref{prop:w2descr}.\square

Abelianization of the equation \Ref{eqn:wipdef} yields
\[\oline{w_j'}=\sum_i\mij\bar w_i.\]
In particular,
\[\oline{w_0'}=\prtl\si=0,\]
so $w_0'$ is a product of commutators, ie
\be{equation*}
w_0'=\prod_{\ell=1}^\rho[u_\ell,u_{\ell+\rho}]
\end{equation*}
for some $u_1,\dots,u_{2\rho}\in F_n$.

We now express $\mu$ in terms of the $u_\ell$'s.
First recall that the oriented model surface $\Si_\rho$ of genus $\rho$ is
obtained from $D^2$ by pairwise identification of certain segments of
$S^1$. These segments parametrize $2\rho$ loops in $\Si_\rho$ that in
turn represent a symplectic basis $b_1,\dots,b_{2\rho}$ for
$H_1(\Si_\rho)$.
Choose a map $f:\Si_\rho\to X$ representing the generator $\si\in
H_2(X)$ such that, in the notation introduced in Section~\ref{sec:formula},
the class $\ti b_\ell\in H_1(X)/T$ is the one
represented by $u_\ell$, for
$\ell=1,\dots,2\rho$. (Compare the map $g$ in the proof of
Proposition~\ref{prop:w2descr}.) Then $\mu$ is given by the formula
\Ref{eqn:mu-formula}.

In the following we use the notation of Section~\ref{sec:hom-int}.
For $i=1,\dots,r-2$ let $z_i\in F_n$ be a lift of $\ga_i\in H$. To prove the
theorem it suffices, by Lemma~\ref{lemma:tphi-eps}, to find a representation
$\psi:F_n\to G=\Sp1$ such that
\[\psi(w'_j)=\eps_j,\quad\psi(z_i)=1\]
for $i=1,\dots,r-2$ and $j=0,\dots,s$. To this end we will apply
Theorem~\ref{thm:b2-1-explicit} with $k=n-2$ and
\[(v_0,\dots,v_{n-2})=(w'_0,\dots,w'_s,z_1,\dots,z_{r-2})\]
and $\eps_j=1$ for $j>s$. Recall that we have fixed a basis for the
free group $\pi_1(X^1,x_0)$, so we can identify $H_1(X^1)=L=\z^n$.
Let $K$ denote the subgroup of $L$ spanned by the linearly independent elements
$\oline{w'_1},\dots,\oline{w'_s}$. By the elementary divisors
theorem we can find a basis $d_1,\dots,d_n$ for $L$ and
integers $m_1,\dots,m_s$ such that $\{m_id_i\}_{1\le i\le s}$ is a
basis for $K$. This implies that $\{d_i\}_{s+1\le i\le n}$ maps to a
basis for $A_0:=H_1(X)/T$.
Therefore,
\be{align*}
\det(\lla\wedge\bar v_1\wedge\dots\wedge\bar v_{n-2})&=
\pm|T|\det(\lla\wedge d_1\wedge\cdots\wedge d_s\wedge
\bar z_1\wedge\dots\wedge\bar z_{r-2})\\
&=\pm|T|\det(\mu\wedge\bar \ga_1\wedge\cdots\wedge\bar \ga_{r-2})\\
&\neq0,
\end{align*}
since the natural map $L\to A_0$ takes $\bar z_i$ to $\bar \ga_i$ and
the induced map $\La^2L\to\La^2A_0$ takes $\lla$ to $\mu$.
Theorem~\ref{thm:b2-1-explicit} now guarantees the existence of the
desired representation $\psi$. This proves 
Theorem~\ref{thm:b2-1}.\square

\appendix

\section{Commutators and group homology}

The purpose of this appendix is to shed some light on the term $\lla$ in
Theorem~\ref{thm:b2-1-explicit}.
Let $G$ be any group and $H_*(G)$ its group homology with integer coefficients.
Let $H_*(G;\q)=H_*(G)\otimes\q$ be
the homology with rational coefficients.
We identify $H_1(G)$ with the abelianization of $G$. The image of an element
$x\in G$ in $H_1(G;\q)$ will be denoted $\check x$. Let $[G,G]$ be the
commutator subgroup of $G$.

\be{prop}
If $H_2(G)$ is a torsion group (or equivalently, if $H_2(G;\q)=0$)
then there exists a group homomorphism
\be{equation}\label{eqn:algg}
\al:[G,G]\to\La^2H_1(G;\q)
\end{equation}
that sends any commutator $[x,y]$ to $\check x\wedge\check y$.
\end{prop}

If $H_2(G)=0$ then the proposition also holds if $H_1(G;\q)$ is replaced
by $H_1(G)$ in \Ref{eqn:algg}.

\proof We use Miller's description \cite{Miller1} (see also
\cite{Ellis1}) of $H_2(G)$ in terms of
commutators. Let $\la G,G\ra$ be the free group on all pairs $\la x,y\ra$
with $x,y\in G$. Let $Z(G)$ be the kernel of the homomorphism
$\la G,G\ra\to[G,G]$ that sends $\la x,y\ra$ to $[x,y]$. Miller shows that
\[H_2(G)\approx Z(G)/B(G),\]
where $B(G)$ is the normal subgroup of $\la G,G\ra$
generated by a certain subset $E$. Now let
\[\beta:\la G,G\ra\to\La^2H_1(G;\q)\]
be the homomorphism that sends $\la x,y\ra$ to $\check x\wedge\check y$.
It is easily verified that $\beta$ vanishes
on $E$, hence on $B(G)$. Since $Z(G)/B(G)$ is a torsion group while
$\La^2H_1(G;\q)$ is torsion free we conclude that $\beta$ vanishes on $Z(G)$.
Hence $\beta$ descends to a homomorphism $\al$ as in the
proposition.\square

Now let $G$ be the free group $F_n$. Then $H_2(G)=0$, so in the notation
of Theorem~\ref{thm:b2-1-explicit} we have
\[\lla=\al(v_0).\]


\begin{thebibliography}{10}

\bibitem{Akbulut-McCarthy}
S.~Akbulut and J.~D. McCarthy.
\newblock {\em Casson's invariant for oriented homology $3$-spheres. An
  exposition}.
\newblock Princeton University Press, 1990.

\bibitem{Baird1}
Th.~J. Baird.
\newblock {Cohomology of the space of commuting $n$-tuples in a compact Lie
  group}.
\newblock {\em Algebr. Geom. Topol.}, 7:737--754, 2007.

\bibitem{Ellis1}
G.~J. Ellis.
\newblock Non-abelian exterior products of groups and exact sequences in the
  homology of groups.
\newblock {\em Glasgow Math. J}, 29:13--19, 1987.

\bibitem{Frohman-Nicas}
C.~Frohman and A.~Nicas.
\newblock An intersection homology invariant for knots in a rational homology
  $3$-sphere.
\newblock {\em Topology}, 33:123--158, 1994.

\bibitem{Gerstenhaber-Rothaus}
M.~Gerstenhaber and O.~S. Rothaus.
\newblock The solution of sets of equations in groups.
\newblock {\em Proc. Natl. Acad. Sci. USA}, 48:1531--1533, 1962.

\bibitem{Gompf-Stipsicz}
R.~E. Gompf and A.~I. Stipsicz.
\newblock {\em $4$-Manifolds and Kirby Calculus}, volume~20 of {\em Graduate
  Studies in Mathematics}.
\newblock American Mathematical Society, 1999.

\bibitem{GHV2}
W.~Greub, S.~Halperin, and R.~Vanstone.
\newblock {\em {Connections, Curvature, and Cohomology}}, volume~II.
\newblock Academic Press, 1973.

\bibitem{Guillemin-Pollack}
V.~Guillemin and A.~Pollack.
\newblock {\em Differential Topology}.
\newblock Prentice-Hall, 1974.

\bibitem{Harper-Saveliev1}
E.~Harper and N.~Saveliev.
\newblock {A Casson-Lin type invariant for links}.
\newblock {\em Pacific J. Math.}, 248:139--154, 2010.

\bibitem{Harper-Saveliev2}
E.~Harper and N.~Saveliev.
\newblock {Instanton Floer homology for two-component links}.
\newblock {\em J. Knot Theory Ramifications}, 21(5), 2012.

\bibitem{Hironaka}
H.~Hironaka.
\newblock Triangulations of algebraic sets.
\newblock {\em Proceedings of Symposia in Pure Mathematics}, 29:165--185, 1975.

\bibitem{Hopf2}
H.~Hopf.
\newblock {\"Uber den Rang geschlossener Liescher Gruppen}.
\newblock {\em Comment. Math. Helv.}, 13:119--143, 1940-41.

\bibitem{Hopf1}
H.~Hopf.
\newblock {\"Uber die Topologie der Gruppen-Mannigfaltigkeiten und ihrer
  Verallgemeinerungen}.
\newblock {\em Annals of Math.}, 42:22--52, 1941.

\bibitem{Lang2}
S.~Lang.
\newblock {\em Algebra}.
\newblock Addison Wesley, third edition, 1993.

\bibitem{X-S-Lin}
Xiao-Song Lin.
\newblock A knot invariant via representation spaces.
\newblock {\em J.~Diff. Geom.}, 35:337--357, 1992.

\bibitem{Miller1}
C.~Miller.
\newblock The second homology group of a group; relations among commutators.
\newblock {\em Proc. Amer. Math. Soc.}, 3:588--595, 1952.

\bibitem{rub-sav0}
D.~Ruberman and N.~Saveliev.
\newblock {Rohlin's invariant and gauge theory, I. Homology $3$-tori}.
\newblock {\em Comment. Math. Helv.}, 79:618--646, 2004.

\bibitem{Spanier}
E.~H. Spanier.
\newblock {\em Algebraic Topology. Corrected reprint}.
\newblock Springer, 1981.

\bibitem{Stallings1}
J.~Stallings.
\newblock Homology and central series of groups.
\newblock {\em J. Algebra}, 2:170--181, 1965.

\bibitem{Thaddeus1}
M.~Thaddeus.
\newblock {A perfect Morse function on the moduli space of flat connections}.
\newblock {\em Topology}, 39:773--787, 2000.

\bibitem{wagon}
S.~Wagon.
\newblock {\em {The Banach-Tarski Paradox}}.
\newblock Cambridge University Press, 1985.

\end{thebibliography}

\noindent\textsc{Centre for Quantum Geometry of Moduli Spaces,\\
Ny Munkegade 118,\\
DK-8000 Aarhus C,\\
Denmark}
\\ \\
Email:\ froyshov@qgm.au.dk

\end{document}